\documentclass[a4paper]{article}

\usepackage{mathptmx}       
\usepackage{helvet}         
\usepackage{courier}        
\usepackage{type1cm}        
\usepackage{makeidx}         
\usepackage{graphicx}        
\usepackage{multicol}        
\usepackage[bottom]{footmisc}

\makeindex             

\usepackage{amsmath}
\usepackage{amssymb}
\usepackage[breaklinks]{hyperref}
\usepackage{subcaption}
\usepackage{float}

\newcommand\blfootnote[1]{%
    \bgroup
    \renewcommand\thefootnote{\fnsymbol{footnote}}%
    \renewcommand\thempfootnote{\fnsymbol{mpfootnote}}%
    \footnotetext[0]{#1}%
    \egroup
}

\date{}

\begin{document}

\title{Next-Gen Gas Network Simulation}
\author{Christian Himpe, Sara Grundel, Peter Benner}
%
%
\maketitle

\abstract{%
To overcome many-query optimization, control, or uncertainty quantification work loads in reliable gas and energy network operations,
model order reduction is the mathematical technology of choice.
To this end, we enhance the model, solver and reductor components of the \texttt{morgen} platform,
introduced in \textsc{Himpe et al} [J.~Math.~Ind. 11:13, 2021],
and conclude with a mathematically, numerically and computationally favorable model-solver-reductor ensemble.}

\blfootnote{Christian Himpe (ORCiD: \href{http://orcid.org/0000-0003-2194-6754}{0000-0003-2194-6754}) \texttt{himpe@mpi-magdeburg.mpg.de} \\
\textbullet~ Sara Grundel (ORCiD: \href{http://orcid.org/0000-0002-0209-6566}{0000-0002-0209-6566}) \texttt{grundel@mpi-magdeburg.mpg.de} \\
\textbullet~ Peter Benner (ORCiD: \href{http://orcid.org/0000-0003-3362-4103}{0000-0003-3362-4103}) \texttt{benner@mpi-magdeburg.mpg.de} \\
Computational Methods in Systems and Control Theory Group at the Max Planck Institute for Dynamics of Complex Technical Systems,
Sandtorstra{\ss}e~1, D-39106 Magdeburg, Germany \\
}

\section{\scalebox{.81}{Model Order Reduction for Gas and Energy Networks}}
Computer-based simulation of gas transport in pipeline networks has been an industrial as well as academic field of interest
since the earliest scientific computing systems~\cite{LotH67}.
Especially, the transient simulation of gas flow and the dynamic gas network behavior are the pinnacle discipline in this regard.
The MATLAB-based \texttt{morgen} -- Model Order Reduction for Gas and Energy Networks -- platform\footnote{See: \url{https://git.io/morgen}}
continues this research by providing a modular open-source software simulation stack for the comparison and benchmarking
of models (discretizations), solvers (time steppers), and reductors (model reduction algorithms) \cite{morHimGB21a}.
Beyond, selecting apposite simulator components or ranking model reduction methods,
an overall goal is the acceleration of forward simulations,
so that many-query tasks relying thereon,
such as optimization, control or uncertainty quantification,
benefit in terms of performance.
In this work, we summarize and enhance the foundational work of \cite{morHimGB21a} with additional details,
and accompany version~\textbf{1.1} of \texttt{morgen}.

\subsection{Modules Overview}
The \texttt{morgen} platform is organized into modules: \emph{models}, \emph{solvers}, \emph{reductors}, \emph{networks} and \emph{tests}.
The \emph{networks} module holds topology and scenario data,
and the \emph{tests} module defines the simulation and model reduction experiments,
thus, we summarize the currently available core modules:
\emph{models}, \emph{solvers}, and \emph{reductors}.
The \emph{models} module assembles a semi-discrete input-output system from a network topology.
Currently, two spatially discrete ODE models are included (Table~\ref{tab:models}).

\pagebreak

\begin{table}[h]\centering
\begin{tabular}{l|l|c|l}
 \textbf{Name} & \textbf{Identifier} & \textbf{port-Hamiltonian?} & \textbf{Reference} \\ \hline
 Midpoint discretization & \texttt{ode\_mid} & No  &\cite[Sec.~2.4.1]{morHimGB21a} \\
 Endpoint discretization & \texttt{ode\_end} & Yes &\cite[Sec.~2.4.2]{morHimGB21a}
\end{tabular}
\caption{Available models in \texttt{morgen} in version~1.1.}
\label{tab:models}
\end{table}

The \emph{solvers} module computes a time-discrete output trajectory from a model and a scenario.
Six ODE solvers are provided in the current version (Table~\ref{tab:solvers}). \\

\begin{table}[h]\centering
\begin{tabular}{l|l|c|l}
 \textbf{Name} & \textbf{Identifier} & \textbf{Comment} & \textbf{Reference} \\ \hline
 Adaptive 2nd Order Rosenbrock  & \texttt{generic} & uses \texttt{ode23s} & \cite[Sec.~5.3.1]{morHimGB21a} \\
 1st Order Implicit-Explicit    & \texttt{imex1}   & non-Runge-Kutta      & \cite[Sec.~5.3.3]{morHimGB21a} \\
 2nd Order Implicit-Explicit    & \texttt{imex2}   & Runge-Kutta          & \cite[Sec.~5.3.4]{morHimGB21a} \\
 Explicit 4th Order Runge-Kutta & \texttt{rk4}     &                      & \cite[Sec.~5.3.2]{morHimGB21a} \\
 Explicit 2nd Order Runge-Kutta & \texttt{rk2hyp}  & increased stability  & \cite{Van72} \\
 Explicit 4th Order Runge-Kutta & \texttt{rh4hyp}  & increased stability  & \cite{MeaR99}
\end{tabular}
\caption{Available solvers in \texttt{morgen} in version~1.1.}
\label{tab:solvers}
\end{table}

The \emph{reductors} module compresses a model given a solver and (generic training) scenario.
All in all, $23$ reductors organized in four classes are available (Table~\ref{tab:reductors}). \\

\begin{table}[h]\centering
\begin{tabular}{l|l|l|l}
 \textbf{Name} & \textbf{Identifier} & \textbf{Linear Variant} & \textbf{Reference} \\ \hline
 Proper Orthogonal Decomposition & \texttt{pod\_r} & -- & \cite[Sec.~4.2]{morHimGB21a} \\
 Empirical Dominant Subspaces & \texttt{eds\_ro} & \texttt{eds\_ro\_l} & \cite[Sec.~4.3]{morHimGB21a} \\
 Empirical Dominant Subspaces & \texttt{eds\_wx} & \texttt{eds\_wx\_l} & \cite[Sec.~4.3]{morHimGB21a} \\
 Empirical Dominant Subspaces & \texttt{eds\_wz} & \texttt{eds\_wz\_l} & \cite[Sec.~4.3]{morHimGB21a} \\ \hline
 Balanced POD & \texttt{bpod\_ro} & \texttt{bpod\_ro\_l} & \cite[Sec.~4.4.3]{morHimGB21a} \\
 Balanced Truncation & \texttt{ebt\_ro} & \texttt{ebt\_ro\_l} & \cite[Sec.~4.4]{morHimGB21a} \\
 Balanced Truncation & \texttt{ebt\_wx} & \texttt{ebt\_wx\_l} & \cite[Sec.~4.4]{morHimGB21a} \\
 Balanced Truncation & \texttt{ebt\_wz} & \texttt{ebt\_wz\_l} & \cite[Sec.~4.4]{morHimGB21a} \\ \hline
 Goal-Oriented POD & \texttt{gopod\_r} & -- & \cite[Sec.~4.5.1]{morHimGB21a} \\
 Balanced Gains & \texttt{ebg\_ro} & \texttt{ebg\_ro\_l} & \cite[Sec.~4.5]{morHimGB21a} \\
 Balanced Gains & \texttt{ebg\_wx} & \texttt{ebg\_wx\_l} & \cite[Sec.~4.5]{morHimGB21a} \\
 Balanced Gains & \texttt{ebg\_wz} & \texttt{ebg\_wz\_l} & \cite[Sec.~4.5]{morHimGB21a} \\ \hline
 DMD Galerkin & \texttt{dmd\_r} & -- & \cite[Sec.~4.6]{morHimGB21a}
\end{tabular}
\caption{Available reductors in \texttt{morgen} in version~1.1.}
\label{tab:reductors}
\end{table}

\section{Enhanced Functionality}
In this section we discuss some properties of the \texttt{morgen} platform.
Specifically, one aspect of each of the core modules (\emph{model}, \emph{solver}, \emph{reductor}) is addressed.
Additionally, further network/scenario data-sets were added in version~1.1, too.

\subsection{Gravity Term}
One component of the gas pipeline model, particularly of the retarding forces in the mass-flux equation, is the gravity term,
which accounts for increase or decrease in momentum due to an incline in a pipeline section.
In \cite{BehBS19} this gravity term is modeled in great detail,
as it does not only consider a height difference between the pipe's end points, as \texttt{morgen} does,
but also the height profile for the full run of the pipe (see \cite[Fig.~11]{BehBS19}).
Both approaches are justified, depending on the aimed accuracy of the model, as discussed in \cite{BacG00}.
Such pipeline height profiles can be included into \texttt{morgen} by supplying a pipe as sequence of virtual pipes,
each connecting two subsequent local height extrema.
Also in \texttt{morgen}~1.1, the gravity term is configurable so it is computable based on the dynamic pressure,
static pressure or not at all.

\subsection{Explicit Solvers}
In \cite{morHimGB21a}, the classic explicit 4th order Runge-Kutta method \texttt{rk4} was tested, as it was employed in earlier works.
Yet, we found it to be not suitable for gas network simulations.
In~\cite{Lew95} an explicit Runge-Kutta method from \cite{Van72} was suggested for this application,
while in \cite{MeaR99} a Runge-Kutta method was optimized in terms of its \linebreak \emph{hyperbolic stability limit}.
The Butcher tableaus for these explicit 5-stage, 2nd order and 6-stage, 4th order methods with increased stability are given by:
\begin{center}
\vskip1em
 \begin{tabular}{c|ccccc}
  $0$ \\
  $\frac{1}{4}$ & $\frac{1}{4}$ \\
  $\frac{1}{6}$ & $0$ & $\frac{1}{6}$ \\
  $\frac{3}{8}$ & $0$ & $0$ & $\frac{3}{8}$ \\
  $\frac{1}{2}$ & $0$ & $0$ & $0$ & $\frac{1}{2}$ \\ \hline
  & $0$ & $0$ & $0$ & $0$ & $1$
 \end{tabular}
\hskip4em
 \begin{tabular}{c|cccccc}
  $0$ \\
  $c_2$ & $a_2$ \\
  $c_3$ & $0$ & $a_3$ \\
  $c_4$ & $0$ & $0$ & $a_4$ \\
  $c_5$ & $0$ & $0$ & $0$ & $a_5$ \\
  $c_6$ & $0$ & $0$ & $0$ & $0$ & $a_6$ \\ \hline
  & $b_1$ & $b_2$ & $b_3$ & $b_4$ & $0$ & $b_6$
 \end{tabular}
\vskip1em
\end{center}

These additional solvers \texttt{rk2hyp}, \texttt{rk4hyp} (see Table~\ref{tab:rkm} for coefficients) were added to \texttt{morgen}~1.1 and tested against various test problems,
and both increased-stability solvers allow larger time-steps then \texttt{rk4},
specifically in conjunction with the \texttt{ode\_end} model,
but compared to the implicit-explicit solvers \texttt{imex1} and \texttt{imex2}, they are still not fully competitive.
However, these explicit methods could be interesting for new implicit-explicit or predictor-corrector methods.

\begin{table}[h]\centering
\vskip1.5em
\begin{tabular}{l|l}
 $c_2 = a_2 = \phantom{-}0.16791846623918$ ~ & ~ $b_1 = -0.15108370762927$ \\
 $c_3 = a_3 = \phantom{-}0.48298439719700$ ~ & ~ $b_2 = \phantom{-}0.75384683913851$ \\
 $c_4 = a_4 = \phantom{-}0.70546072965982$ ~ & ~ $b_3 = -0.36016595357907$ \\
 $c_5 = a_5 = \phantom{-}0.09295870406537$ ~ & ~ $b_4 = \phantom{-}0.52696773139913$ \\
 $c_6 = a_6 = \phantom{-}0.76210081248836$ ~ & ~ $b_6 = \phantom{-}0.23043509067071$
\end{tabular}
\caption{Butcher tableau coefficients for the \texttt{rk4hyp} method; values taken \linebreak from~\cite[Sec.~4.2]{MeaR99}.}
\label{tab:rkm}
\end{table}

\pagebreak

\subsection{Gain Matching}
An important quality for certain applications of model reduction,
such as electrical circuits, is the preservation of the steady-state gain (also known as DC gain), which is the output for zero frequency input.
First, we clarify that we are not discussing the actual steady-state gain of the reduced order model,
due to the centering around the steady-state and hence, the steady-state gain match \cite[Sec.~3]{morHimGB21a}.
Yet, there can still be an output error for a constant input on top of the steady-state input,
which is relevant due to the assumed low-frequency boundary values.
Since there is an interpretation of gas networks as circuits \cite{TaoT98}, we consider this reduced model property,
which induces two questions:
How to compute the steady-state gain, and how to correct a gain mismatch?
The former is answered by \cite{van16},
stating that for a linear port-Hamiltonian model, with components as in \cite[Sec.~2.9]{morHimGB21a}, the gain is computable~by:
\begin{align*}
 S = C Q^{-1} B.
\end{align*}

Since the models are nonlinear and do not have to be port-Hamiltonian, but comprise the same model components,
the above formula can still be applied albeit yielding only an approximation. 
The per-port gain mismatch is then computed by the difference of full and reduced order model gain:
\begin{align*}
 D := (C Q^{-1} B) - (C_r Q_r^{-1} B_r),
\end{align*}
which can then be used to correct the reduced order model gain by adding it as a feedthrough matrix to the output function,
as described in the gain matching procedure in \cite{morSamPG95}.
We added this approximate gain matching test to \texttt{morgen}~1.1.

The gain correction was tested with all reductors (Table~\ref{tab:reductors}).
For all reductors the correction was about the level of $10^{-5}$,
except for the \texttt{bpod\_ro} method, for which the gain correction fully deteriorates the reduced order model.
Thus, the improvement of reduced order models is small at best.
This is not unexpected, considering the gas network model is hyperbolic:
A single pipeline, or more generally an input-output system based on a first order hyperbolic partial differential equation,
has the transport property which expresses as a delay in observable outputs of controllable inputs.
Hence, an immediate transformation of inputs to outputs (circumventing the system dynamics),
as a feedthrough term does, is typically not needed.

\section{Numerical Experiments}
\begin{figure}[H]

\begin{subfigure}{0.49\textwidth}
\includegraphics[width=\textwidth]{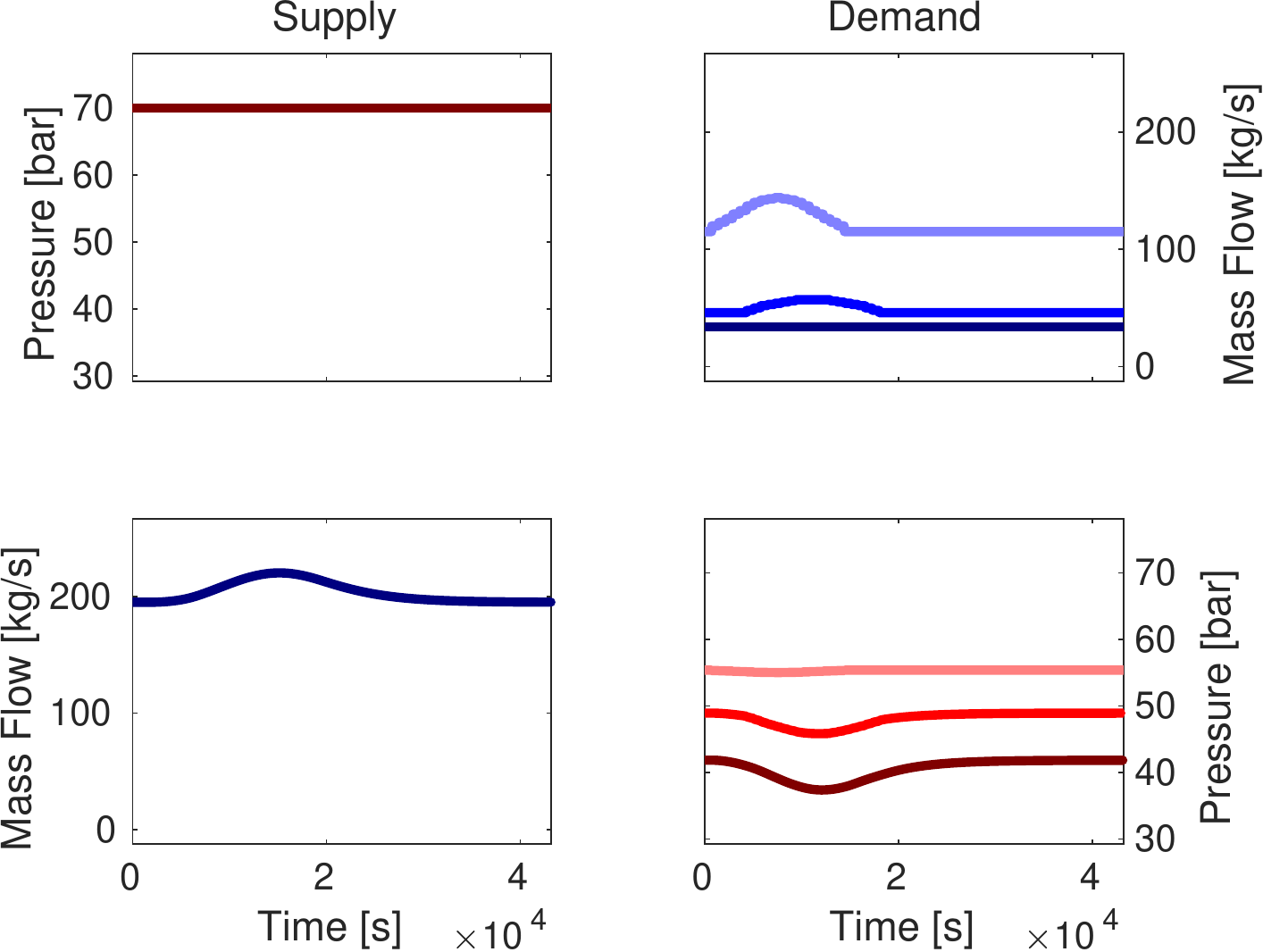}
\caption{\textbf{Hypothetical} network's test scenario.}
\label{fig:numex1:scen}
\end{subfigure}
~~
\begin{subfigure}{0.48\textwidth}
\includegraphics[width=\textwidth]{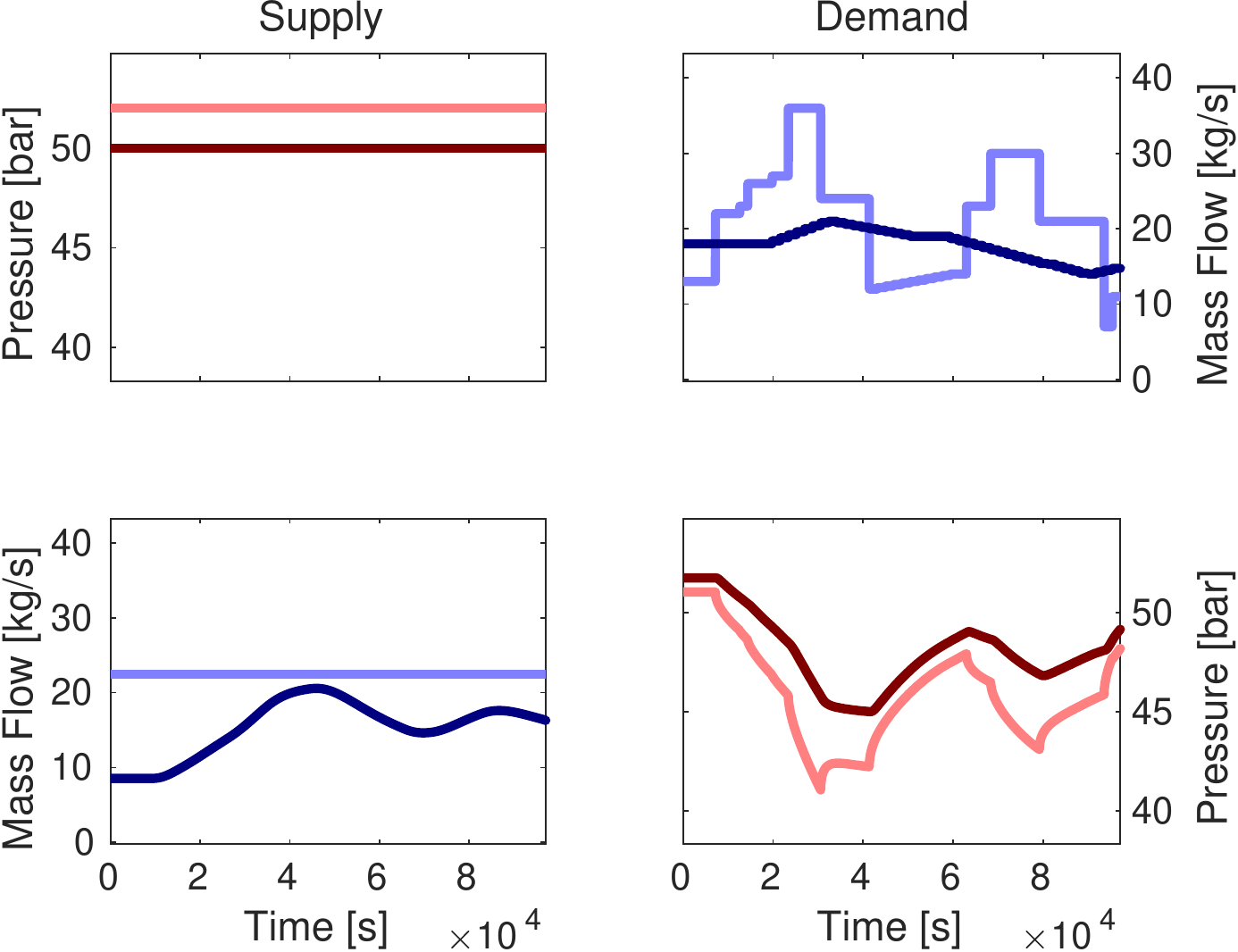}
\caption{\textbf{Actual} network's test scenario.}
\label{fig:numex2:scen}
\end{subfigure}

\begin{subfigure}{0.49\textwidth}
\includegraphics[width=\textwidth]{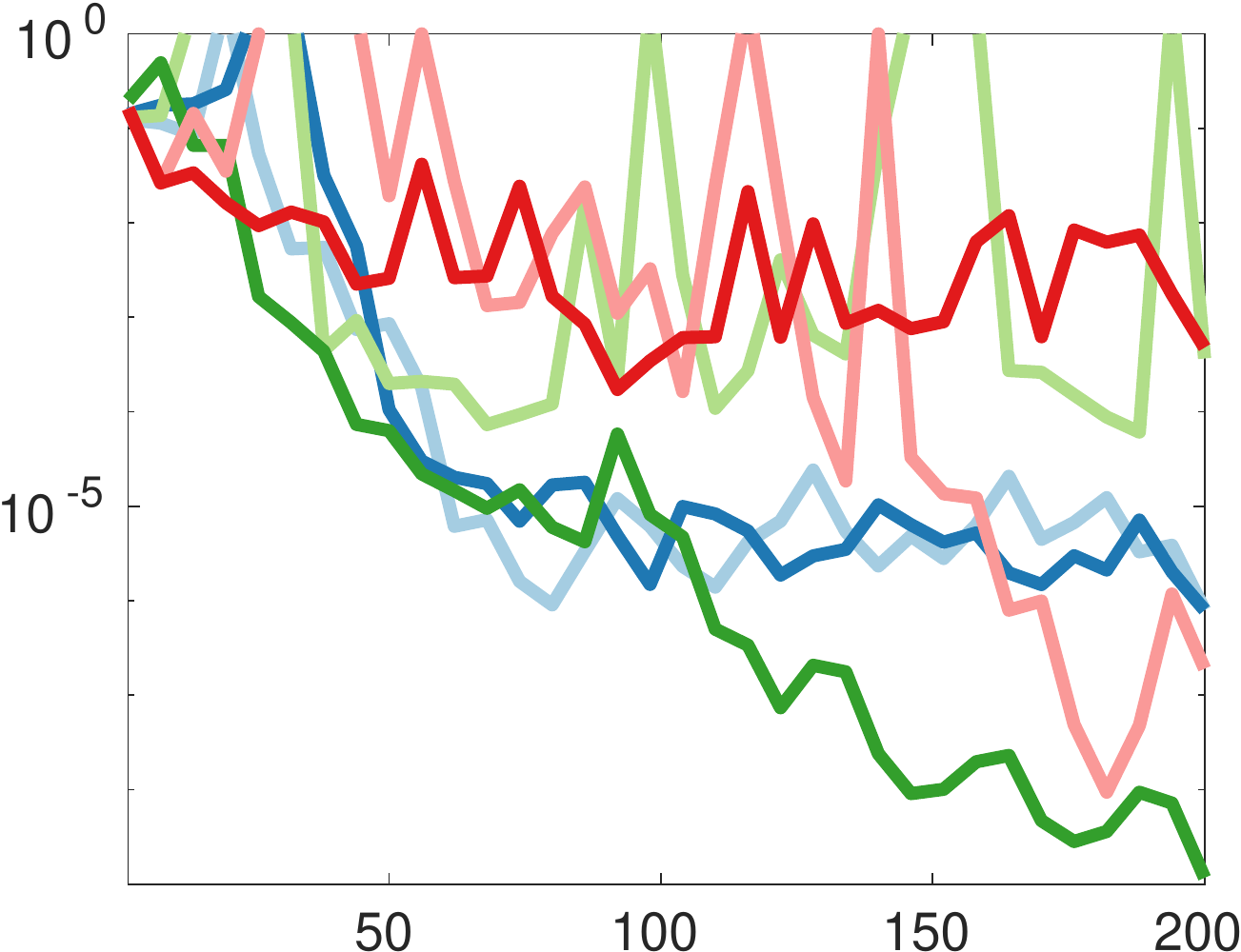}
\caption{Relative $L_2 \otimes L_2$ error between ROM and FOM for the \texttt{ode\_end} model, \texttt{imex1} solver, and \emph{linear} reductors versus reduced order for the \textbf{hypothetical} network.}
\label{fig:numex1:end1}
\end{subfigure}
~~
\begin{subfigure}{0.49\textwidth}
\includegraphics[width=\textwidth]{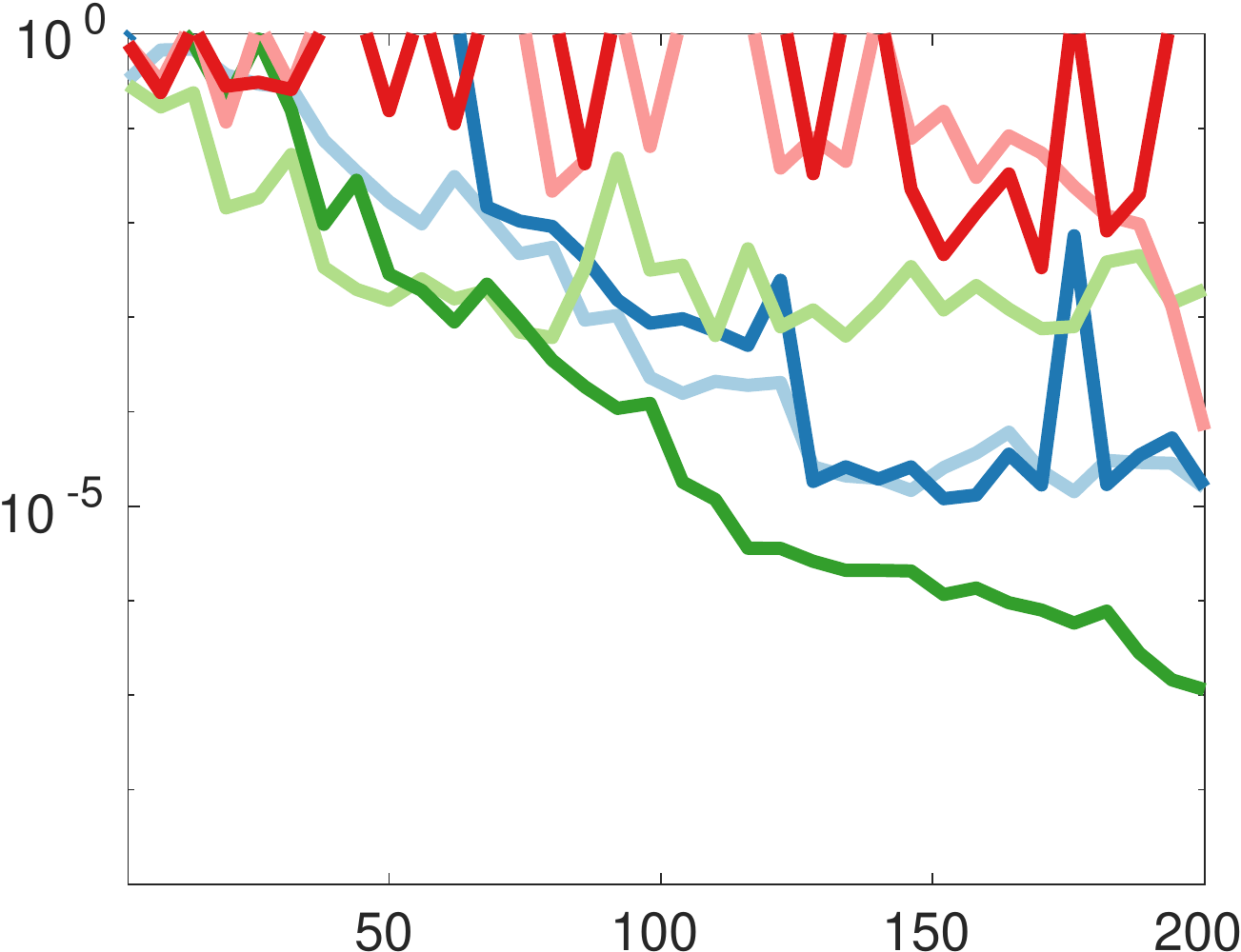}
\caption{Relative $L_2\,\otimes\,L_2$ error between ROM and FOM for the \texttt{ode\_end} model, \texttt{imex1} solver, and \emph{linear} reductors versus reduced order for the \textbf{actual} network.}
\label{fig:numex2:end1}
\end{subfigure}

~\\

\begin{subfigure}{\textwidth}\centering
\includegraphics[width=.5\textwidth]{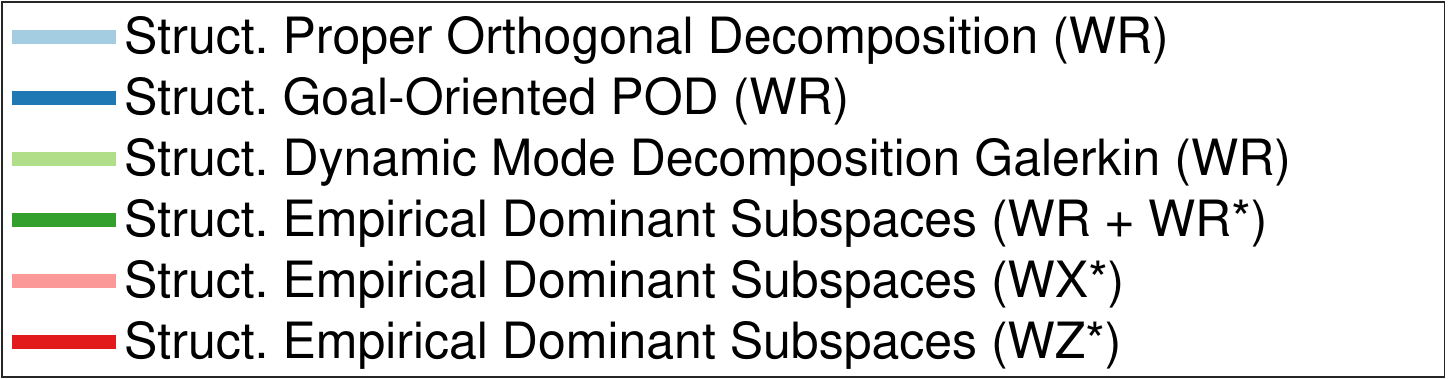}
\caption{Common legend for the model reduction error plots.}
\label{fig:numex1:legend}
\end{subfigure}

~\\

\begin{subfigure}{0.48\textwidth}\footnotesize
\begin{tabular}{r|c|c}
 \textbf{Reductor} & \textsc{MORscore} & \textbf{Avg. Gain Error} \\ \hline
 \texttt{pod\_r}     & $0.27$ & $6 \cdot 10^{-6}$ \\
 \texttt{gopod\_r}   & $0.26$ & $6 \cdot 10^{-6}$ \\
 \texttt{dmd\_r}     & $0.18$ & $8 \cdot 10^{-6}$ \\
 \texttt{eds\_ro\_l} & $0.30$ & $8 \cdot 10^{-6}$ \\
 \texttt{eds\_wx\_l} & $0.18$ & $8 \cdot 10^{-6}$ \\
 \texttt{eds\_wz\_l} & $0.15$ & $8 \cdot 10^{-6}$ \\
\end{tabular}
\caption{\textsc{MORscore}s $\mu(200,\epsilon_{\operatorname{mach}(16)})$ in the \linebreak $L_2 \otimes L_2$ error norm,
         and mean steady-state gain error for the \textbf{hypothetical} network.}
\label{fig:numex1:mid2}
\end{subfigure}
~~
\begin{subfigure}{0.48\textwidth}\footnotesize
\begin{tabular}{r|c|c}
 \textbf{Reductor} & \textsc{MORscore} & \textbf{Avg. Gain Error} \\ \hline
 \texttt{pod\_r}     & $0.19$ & $2 \cdot 10^{-5}$ \\
 \texttt{gopod\_r}   & $0.15$ & $1 \cdot 10^{-5}$ \\
 \texttt{dmd\_r}     & $0.15$ & $2 \cdot 10^{-5}$ \\
 \texttt{eds\_ro\_l} & $0.24$ & $2 \cdot 10^{-5}$ \\
 \texttt{eds\_wx\_l} & $0.04$ & $2 \cdot 10^{-5}$ \\
 \texttt{eds\_wz\_l} & $0.03$ & $2 \cdot 10^{-5}$ \\
\end{tabular}
\caption{\textsc{MORscore}s $\mu(200,\epsilon_{\operatorname{mach}(16)})$ in the \linebreak $L_2 \otimes L_2$ error norm,
         and mean steady-state gain error for the \textbf{actual} network.}
\label{fig:numex1:end2}
\end{subfigure}

\caption{Visualization of the test scenario, model reduction errors, \textsc{MORscore}s, and gain errors of the tested ROMs for the
\textbf{hypothetical} network \cite[Part~2]{LotH67} (left side) and \textbf{actual} network \cite[Part~3]{LotH67} (right side).
Computed with MATLAB 2021a.}
\label{fig:numex}
\end{figure}

We extend the numerical experiments in \cite{morHimGB21a}, by reimplementing the results from~\cite{LotH67},
specifically we test the \textbf{hypothetical} network \cite[Part~2]{LotH67},
and the \textbf{actual} network \cite[Part~3]{LotH67} on their associated scenarios.
Both are tree networks, and the empirical-Gramian-based Galerkin reductors \texttt{pod\_r}, \texttt{gopod\_r}, \texttt{dmd\_r},
\texttt{eds\_ro\_l}, \linebreak \texttt{eds\_wx\_l}, \texttt{eds\_wz\_l} are tested on the port-Hamiltonian endpoint model \linebreak \texttt{ode\_end} and the first order implicit-explicit solver \texttt{imex1}.
The results are presented in Figure~\ref{fig:numex}.
In line with other experiments, the \texttt{eds\_ro\_l} reductor yields the most accurate results.

\section{Next-Gen Gas Network Simulation}
Based on the heuristic comparison in \cite{morHimGB21a} and this work's numerical results,
we recommend a port-Hamiltonian model, an implicit-explicit solver, and a Galerkin reductor,
particularly, the endpoint discretization, the first order IMEX time stepper, and the structured empirical dominant subspaces method
as the model-solver-reductor ensemble, for the next generation of transient gas network simulators.

\bibliographystyle{plainurl}
\bibliography{mor,csc,software}

\begin{thebibliography}{10}

\bibitem{BacG00}
G.~Bachman and M.~Goodreau.
\newblock Less is more accuracy versus precision in modeling.
\newblock In {\em PSIG Annual Meeting 2000}, pages PSIG--0009, 2000.
\newblock URL:
  \url{https://onepetro.org/PSIGAM/proceedings-abstract/PSIG00/All-PSIG00/PSIG-0009/2043}.

\bibitem{BehBS19}
M.~Behbahani-Nejad, A.~Berm{\'u}dez, and M.~Shabani.
\newblock Finite element solution of a new formulation for gas flow in a pipe
  with source terms.
\newblock {\em Journal of Natural Gas Science and Engineering}, 61:237--250,
  2019.
\newblock \href {https://doi.org/10.1016/j.jngse.2018.11.019}
  {\path{doi:10.1016/j.jngse.2018.11.019}}.

\bibitem{morHimGB21a}
C.~Himpe, S.~Grundel, and P.~Benner.
\newblock Model order reduction for gas and energy networks.
\newblock {\em Journal of Mathematics in Industry}, 11:13, 2021.
\newblock \href {https://doi.org/10.1186/s13362-021-00109-4}
  {\path{doi:10.1186/s13362-021-00109-4}}.

\bibitem{Lew95}
A.~Lewandowski.
\newblock New numerical methods for transient modeling of gas pipeline
  networks.
\newblock In {\em PSIG Annual Meeting}, pages PSIG--9510, 1995.
\newblock URL:
  \url{https://onepetro.org/PSIGAM/proceedings-abstract/PSIG95/All-PSIG95/PSIG-9510/2571}.

\bibitem{LotH67}
L.~A. Lotito and P.~W. Halbert.
\newblock Computer simulation of gas flow dynamics.
\newblock {\em Pipeline Engineer}, 39:31--33, 29--31, 45--47, 1967.

\bibitem{MeaR99}
J.~L. Mead and R.~A. Renaut.
\newblock Optimal {R}unge-{K}utta methods for first order pseudospectral
  operators.
\newblock {\em Journal of Computational Physics}, 152(1):404--419, 1999.
\newblock \href {https://doi.org/10.1006/jcph.1999.6260}
  {\path{doi:10.1006/jcph.1999.6260}}.

\bibitem{morSamPG95}
R.~Samar, I.~Postlewaite, and D.~W. Gu.
\newblock Model reduction with balanced realizations.
\newblock {\em Internat. J. Control}, 62(1):33--64, 1995.
\newblock \href {https://doi.org/10.1080/00207179508921533}
  {\path{doi:10.1080/00207179508921533}}.

\bibitem{TaoT98}
W.~Q. Tao and H.~C. Ti.
\newblock Transient analysis of gas pipeline network.
\newblock {\em Chemical Engineering Journal}, 69(1):47--52, 1998.
\newblock \href {https://doi.org/10.1016/S1385-8947(97)00109-5}
  {\path{doi:10.1016/S1385-8947(97)00109-5}}.

\bibitem{Van72}
P.~J. van~der Houwen.
\newblock Explicit {R}unge-{K}utta formulas with increased stability
  boundaries.
\newblock {\em Numerische Mathematik}, 20:149--164, 1972.
\newblock \href {https://doi.org/10.1007/BF01404404}
  {\path{doi:10.1007/BF01404404}}.

\bibitem{van16}
A.~van~der Schaft.
\newblock Interconnections of input-output {H}amiltonian systems with
  dissipation.
\newblock In {\em Proceedings of the 55th IEEE Conference on Decision and
  Control}, pages 4886--4691, 2016.
\newblock \href {https://doi.org/10.1109/CDC.2016.7798983}
  {\path{doi:10.1109/CDC.2016.7798983}}.

\end{thebibliography}
\end{document}